\documentclass[12pt,leqno]{article}

\usepackage{amsmath,amssymb}
\usepackage{latexsym}
\usepackage{graphicx}
\usepackage{bm}

\def\R{\mathbb{R}}
\def\C{\mathbb{C}}

\begin{document}

\baselineskip=16pt

\title{A solution to the Pompeiu problem}

\author{ A. G. Ramm\\
\\
Mathematics Department, Kansas State University,\\
Manhattan, KS 66506-2602, USA\\
email: ramm@math.ksu.edu }

\renewcommand{\thefootnote}{\fnsymbol{footnote}}

\date{}
\maketitle
\begin{abstract}
\baselineskip=14pt
Let $f \in L_{loc}^1 (\R^n)\cap \mathcal{S}$, where
$\mathcal{S}$ is the Schwartz class of distributions,  and
$$\int_{\sigma (D)} f(x) dx = 0 \quad \forall \sigma \in G, \qquad (*)$$
where  $D\subset \R^n$ is a bounded domain,  the closure $\bar{D}$ of
which is diffeomorphic to a closed ball, and $S$ is its boundary. Then the complement of $\bar{D}$
is connected and path connected.
By $G$ the group of all rigid motions of $\R^n$ is denoted. This group
consists of all translations and rotations.
A proof of the following theorem is given.

Theorem 1. {\it Assume that $n=2$, $f\not\equiv 0$,  and (*)
holds. Then $D$ is a ball.}

Corollary.  {\it If the problem $(\nabla^2+k^2)u=0$ in $D$,
$u_N|_S=0$, $u|_S=const\neq 0$ has a solution, then $D$ is a ball.}

Here $N$ is the outer unit normal to $S$. 
\end{abstract}



\section{Introduction}
Let $f \in L_{loc}^1 (\R^n)\cap \mathcal {S}$, where
$\mathcal{S}$ is the Schwartz class of distributions,  and
\begin{equation}
\label{e1}
\int_{\sigma (D)} f(x) dx = 0 \quad \forall \sigma \in G,
\end{equation}
where $G$  is the group of all rigid motions of $\R^n$, $G$  consists of
all translations and rotations, and
 $D\subset \R^n$ is a bounded domain, the closure $\bar{D}$ of
which is diffeomorphic to a closed ball. Under these assumptions the
complement of $\bar{D}$
in $\R^n$ is connected and path connected (\cite{H}).
By $S$ the boundary of $D$ is denoted, and  $N$ denotes the unit normal 
to $S$ pointing out of $D$.
In \cite{P} the following question was raised by D.Pompeiu:

{\it Does (1) imply that  $f = 0$?}

If yes, then we say that  $D$  has $P$-property (Pompeiu's
property), and write $D\in P$. Otherwise, we say
that $D$ fails to have $P$-property, and write  $D\in \overline{P}$.
Pompeiu claimed that every plane bounded domain has $P-$property, but
a counterexample was given 15 years later in \cite{C}. The counterexample
is a
domain $D$ which is a disc, a ball in $\R^n$ for $n> 2$. If $D$ is a
ball, then there are $f\not\equiv 0$ for which equation (1) holds.
The set of all $f\not\equiv 0$, for which equation
(1) holds, was constructed in \cite{R363}. A bibliography on the Pompeiu
problem ($P-$problem) can be found in \cite{Z}. The results
on $P-$problem which are used in this paper are derived
in \cite{R635}. The $P-$problem is equivalent to a symmetry problem,
see Corollaries 1,2 below.
The author's results on other symmetry
problems are given in \cite{R512} and \cite{R626}.
The modern formulation of the $P-$problem is the following:

{\it Prove that if $D\in \overline{P}$
then $D$ is a ball.}

We use the word ball also in the case $n=2$, when this word means
disc, and solve the $P-$problem. The proof of Theorem 1 we give assuming $n=2$,
but this proof is easily generalized to the case $n>2$.
Our standing assumptions are:

{\it Assumptions A:
a) $D$ is a bounded domain, the closure of which is diffeomorphic to a
closed ball, the boundary $S$ of $D$ is a closed connected $C^1-$smooth
surface, b) $D$ fails to have $P-$property, and c) $n=2$.}

{\bf Theorem 1.} {\it If Assumptions A hold, then $D$ is a ball.}

{\bf Corollary 1.} {\it If problem \eqref{e3} (see below) has a solution,
then $D$ is a ball.}

{\bf Corollary 2.} {\it If the problem $(\nabla^2+k^2)u=0$ in $D$,
$u_N|_S=0$, $u|_S=const\neq 0$ has a solution, then $D$ is a ball.}

In Section 2 these results are proved.

\section{Proof of Theorem 1}

If Assumptions A hold, then the boundary $S$ of $D$
is real-analytic (see \cite{W}) and 
\begin{equation}
\label{e2}
\int_De^{ik\alpha\cdot x}dx=0,\quad \forall \alpha\in S^1,
\end{equation}
where $S^1$ is the unit sphere in $\R^2$, and $k>0$ is a fixed number, see 
\cite{R635}.

The following Lemmas 1-3 are proved in \cite{R635}
(Lemma 1 is Lemma 3 in \cite{R635}, Lemma 2 is Lemma 5 in \cite{R635},
and Lemma 3 is formula  (32) in \cite{R635}):

{\bf Lemma 1.} {\it If and only if relation  
\eqref{e2} holds then the overdetermined problem
\begin{equation}
\label{e3}
(\nabla^2+k^2)u=1 \quad in \quad D, \quad u|_S=0, \quad u_N|_{S}=0,
\end{equation}
has a solution. }

{\bf Lemma 2.} {\it If \eqref{e2} holds for all $\alpha\in S^1$ then it
holds for all $\alpha\in M$, where $M:=\{z: z\in\C^2, z_1^2+z_2^2=1\}$.}

The $M$ is an algebraic variety intersecting $\R^2$ over $S^1$.

Let us assume that the boundary $S$ is star-shaped.
Let $r=f(\phi)$ be the equation of $S$, 
where $0<c_2\le f\le c_2$,  $c_j$ are constants, $j=1,2,$
and $f$ is a smooth $2\pi-$periodic function.

{\bf Lemma 3.} {\it If \eqref{e2} holds for all $\alpha\in S^1$, then
\begin{equation}\label{e4}
\int_{-\pi}^{\pi} f'(\phi)f(\phi)e^{ikf(\phi)\cos(\phi-\theta)}d\phi=0,
\quad \forall \theta\in \C.
\end{equation}
}
Let us choose $\cos\theta=is$ and $\sin\theta=(s^2+1)^{1/2}$. Then
$\{is, (s^2+1)^{1/2}\}\in M$, and \eqref{e4} can be written as
\begin{equation}\label{e5}
\int_{-\pi}^{\pi}
f'(\phi)f(\phi)e^{-skf(\phi)\cos\phi+ik(s^2+1)^{1/2}f(\phi)\sin\phi}d\phi=0,
\quad \forall s>0.
\end{equation}
Multiply \eqref{e5} by $e^{-As}$, where $A>0$ is a large constant, and
integrate over $s$ from $0$ to $\infty$. Then one gets
\begin{equation}\label{e6}
\int_{-\pi}^{\pi}d\phi f'(\phi)f(\phi)\int_0^\infty ds
e^{-s(a+A)+i(s^2+1)^{1/2}b}=0, \quad \forall A>A_0,
\end{equation}
where $A_0>0$ is a fixed large constant, 
$$a=a(\phi)=kf(\phi)\cos\phi, \quad
b=b(\phi)=k f(\phi)\sin\phi, \quad A_0>\max_{\phi\in [-\pi,\pi]}|a(\phi)|.$$

One has
\begin{equation}\label{e7}
\int_0^\infty e^{-s(a+A)+i(s^2+1)^{1/2}b}ds=
(a+A)^{-1}e^{ib}[1+O(A^{-1})]ds,\quad A\to \infty.
\end{equation}
Writing $$(a+A)^{-1}=\sum_{j=0}^\infty(-1)^ja^jA^{-1-j},\quad A>A_0,$$
one obtains from \eqref{e6} and \eqref{e7} the relation
\begin{equation}\label{e8}
\int_{-\pi}^{\pi}
f'(\phi)f(\phi)e^{ib}\sum_{j=0}^\infty(-1)^ja^jA^{-1-j}[1+O(A^{-1})]
d\phi=0,\quad  A\to \infty.
\end{equation}
Multiply \eqref{e8} by $A$ and let $A\to \infty$. This yields relation
\eqref{e9}, see below, with $j=0$. After getting relation 
\eqref{e9} with $j=0$, multiply \eqref{e8} by $A^2$ and let $A\to \infty$.
This yields relation \eqref{e9} with $j=1$.
Continue in 
this fashion to get
\begin{equation}\label{e9}
\int_{-\pi}^{\pi} f'(\phi)f(\phi)a^je^{ib}d\phi=0,\quad\forall j=0,1,....
\end{equation}
Applying the Laplace method (see \cite{Fe}) for calculating the asymptotic
behavior of integral \eqref{e9} as $j\to \infty$, one
concludes that \eqref{e9} can hold if and only if $f'=0$,
that is, if and only if $f=const$.

Let us give details. Consider the function $a^{2m}=e^{m\Psi}$, where 
$\Psi:=\ln[k^2f^2(\phi)\cos^2\phi]$, $j=2m$, so that the expression under
the logarithm sign is non-negative.
The stationary points of the function $\Psi$ are found from
the equation $\frac{f'(\phi)}{f(\phi)} -\tan\phi=0$. 

If $D$ is not a ball, then the function $f(\phi)$ attains its maximum $F$
at a point, which one may denote $\phi=0$. There can be finitely many
points at which $f$ attains local maximums, because $f$ is analytic.
There are finitely many points at which $f$ attains the value $F$.
We assume for simplicity that these points are non-degenerate, so
$f^{''}< 0$ at these points. Since $f>0$, one has  
the inequality 
$$\frac{d}{d\phi}\Big(\frac{f'(\phi)}{f(\phi)}-\tan\phi \Big)=
\frac{f^{''}(\phi)}{f(\phi)}-\frac{(f'(\phi))^2}{f^2(\phi)}-
\frac{1}{\cos^2\phi}<0,$$
if $f^{''}< 0$. Therefore, the critical points are non-degenerate and 
 the main term of the asymptotic of the integral
\eqref{e9} with $j=2m$ as $m\to \infty$, corresponding to
the stationary point $\phi=0$ can be calculated as follows.
Let $I$ denote the integral in \eqref{e9}.
The stationary point $\phi=0$ is a non-degenerate interior
point of maximum of $f$ and, therefore, of $\Psi$. 
Since $e^{ib(\phi)}=1+ ikf(\phi)\sin\phi+...$,
$f(\phi)=f(0)+f'(0)\phi+f^{''}(0)\phi^2/2+...$ and
$f'(\phi)=f'(0)+f^{''}(0)\phi+f^{'''}(0)\phi^2/2+...$,
$\Psi(\phi)=\Psi(0)-\gamma \phi^2+....$, where $\gamma:=|\Psi^{''}(0)|$, 
one multiplies the three terms $ff'e^{ib}$, takes into account that 
$f'(0)=0$ and $\Psi'(0)=0$ at the critical point, and gets  
$I\sim e^{m\Psi(0)}J$, where
$$J=\int_{[-\delta, 
\delta]}\Big((ikf^2(0)f^{''}(0)+f(0)f^{'''}(0)/2 
)\phi^2+...\Big)e^{-m\gamma \phi^2}d\phi.$$ 
As $m\to \infty$, one extends the interval of integration to $(-\infty, 
\infty)$ and calculates the main term of the asymptotic of $J$ as
$m\to \infty$ by using the formula 
$$\int_{-\infty}^\infty\phi^2e^{-m\gamma \phi^2}d\phi=\frac{\Gamma(3/2)}
{(m\gamma)^{3/2}},$$
 where  $\Gamma(z)$ is the Gamma-function, $\Gamma(3/2)=\sqrt{\pi}/2$.
The result is
\begin{equation}\label{e10}
I\sim 
e^{m\Psi(0)}\frac{\Gamma(3/2)}{(m\gamma)^{3/2}}\Big(ikf^2(0)f^{''}(0)+
f(0)f^{'''}(0)/2\Big),\quad m\to \infty. 
\end{equation}
Since $I=0$ and $f(0)>0$, one concludes from \eqref{e10}, after taking the 
imaginary part, that 
$f^{''}(0)=0$, and after taking the real part, that $f^{'''}(0)=0$. This 
contradicts the non-degeneracy of the critical point 
$\phi=0$. If one does not assume the non-degeneracy of this critical 
point, then one uses the analyticity of the function $f$ and concludes 
that if for some $j$ the derivative $f^{(j)}(0)\neq 0$, then  this leads
to a contradiction. Thus, all the derivatives $f^{(j)}(0)=0$ for $j>0$.
Each critical point at which $f=F$ can be taken to be the point
$\phi=0$, because  the origin for $\phi$ in formula  \eqref{e4}
can be chosen arbitrarily on the interval of length of the period 
$[0,2\pi]$.

If the critical point $\phi=0$ is non-degenerate then 
the inputs of local maximums at which $f=F$ cannot compensate each
other since their imaginary parts are all of the same sign 
since $f^{''}< 0$  and $f>0$ at these points. There can be at 
most finitely many critical points of $f$ since $f$ is analytic.  

Thus, the only possibility to have equalities \eqref{e9} for all
large $j$ is to have $f=const$.

Theorem 1 is proved.\hfill $\Box$


\end{document}